\documentclass[12pt]{article}
\usepackage{graphicx}
\usepackage{amsmath, amssymb}

\newtheorem{theorem}{Theorem}
\newtheorem{lemma}[theorem]{Lemma}
\newtheorem{prop}[theorem]{Proposition}

\newtheorem{df}[theorem]{Definition}

\title{Expansive periodic mechanisms}
\author{Ciprian S. Borcea and Ileana Streinu}
\date{}

\begin{document}
\maketitle

\begin{abstract}
A one-parameter deformation of a periodic bar-and-joint framework is expansive when all distances between joints increase or stay the same. In dimension two, expansive behavior can be fully explained through our theory of periodic pseudo-triangulations. However, higher
dimensions present new challenges. In this paper we study a number of periodic
frameworks with expansive capabilities in dimension $d\geq 3$ and register both similarities and contrasts with the two-dimensional case.
\end{abstract}

\medskip

\section{Introduction}

In this paper we investigate expansive capabilities in periodic bar-and-joint frameworks.
The expansive property  refers to one-parameter deformations, sometimes simply called mechanisms, which, for any increase of the parameter, increase or maintain the distance
between any pair of joints. By definition, a bar-and-joint framework will have expansive
capabilities when possessing some non-trivial one-parameter expansive deformation.

\medskip
In dimension two, expansive motions have important applications in robot arm motion planning 
and linkage problems, \cite{S1},\cite{S2}. Mathematically, the key underlying structure is that of
{\em pointed pseudo-triangulation}, \cite{S2} and properties of planar pseudo-triangulations
have been widely explored, \cite{RSS1}, \cite{RSS2}. As a result of our periodic version of
Maxwell's lifting theorem, \cite{M}, \cite{BS4}, expansive planar periodic mechanisms can be similarly
understood in terms of periodic pseudo-triangulations,
 \cite{BS4},\cite{BS6}.

\medskip
However, in dimensions higher than two, a structural understanding of expansiveness remains,
in many respects, an open problem. In the present study, we establish some necessary conditions
and explore a number of periodic families of frameworks with expansive capabilities in dimension three and higher. 

\medskip
After a brief review of basic notions regarding periodic bar-and-joint frameworks and their
deformations \cite{BS1}, \cite{BS2}, we find in Section~\ref{analog} the proper analog for the planar pointedness 
condition and refer to it as {\em pointedness in codimension two}. Guided by this necessary condition, we define and investigate in Section~\ref{type} a type of periodic framework which can be constructed in arbitrary dimension and has expansive capabilities.  We determine the {\em cone of infinitesimally expansive deformations} for this family. In Section~\ref{stress} we
illustrate some new features of expansive behavior in dimension higher than two by investigating a stressed structure. Section~\ref{end} recalls the fact that expansive implies auxetic and presents our conclusions.

\medskip

\section{Periodic frameworks and deformations}\label{preamble}

\vspace{0.2in}
We recall some fundamental definitions on periodic frameworks and their deformation theory, as 
introduced in \cite{BS1} and developed in \cite{BS2}, \cite{BS3}.

A $d$-periodic framework, denoted as $(G,\Gamma, p,\pi)$, is given by an infinite graph $G$, a periodicity group $\Gamma$ acting on $G$, a placement $p$ and a representation $\pi$. The graph $G=(V,E)$ is simple (has no multiple edges and no loops) and connected, with an infinite set of vertices $V$  and (unoriented) edges $E$. The {\em periodicity group} $\Gamma\subset Aut(G)$ is a free Abelian group of rank $d$, acting on $G$ without fixed vertices or fixed edges.  We assume that the quotient  multigraph $G/\Gamma$ (which may have loops and multiple edges) is {\em finite}, and use $n=|V/\Gamma|$  and $m=|E/\Gamma|$ to denote the number of vertex and edge orbits.  The function $p:V\rightarrow R^d$ gives a specific placement of the vertices as points in the Euclidean space $R^d$, in such a way that any two vertices joined by an edge in $E$ are mapped to distinct points. The {\em injective group morphism} $\pi: \Gamma\rightarrow  {\cal T}(R^d)$ gives a faithful representation of $\Gamma$ by a {\em lattice of translations} $\pi(\Gamma)=\Lambda$ of rank $d$ in the group of translations ${\cal T}(R^d)\equiv R^d$.  The placement is {\em periodic} in the obvious sense that the abstract action of the periodicity group $\Gamma$ is replicated by the action of the periodicity lattice $\Lambda=\pi(\Gamma)$
on the placed vertices: $p(\gamma v)=\pi(\gamma) (p(v)), \ \ \mbox{for all}\ \gamma\in \Gamma, v\in V$.

\vspace{0.2in}
A one-parameter {\em periodic deformation} is a family of placements and a family of lattices parametrized by time $(p(t), \pi(t))_t$, such that all bar lengths are maintained and the same abstract periodicity group $\Gamma$  acts on all the frameworks of the deformed family. A periodic framework is {\em rigid} if it has no periodic deformations other than the trivial ones resulting from Euclidean isometries. The configuration space of the periodic framework is obtained from the placements of vertex orbits, subject to the algebraic constraints of prescribed (squared) lengths for edges. We factor out the action of the ${d+1\choose 2}$-dimensional group of Euclidean isometries.  A framework is rigid when corresponding to an isolated point of the configuration space; otherwise it is flexible. The deformation space of a framework is the connected component of the configuration space to which it belongs. 

\vspace{0.2in}
After choosing vertex representatives for all vertex orbits and $d$ generators for the periodicity lattice, the Jacobian matrix at a given placement $p$ is a $(dn+d^2)\times m$ matrix (called the {\em periodic rigidity matrix}) whose rank thus cannot exceed $dn+{d\choose 2}$. At a regular point the rank of the Jacobian equals the dimension of the configuration space in a small neighborhood. We say that a periodic framework is {\em infinitesimally rigid} if its periodic rigidity matrix has the maximum rank of $dn+{d\choose 2}$. In this case, the framework must have at least $dn+{d\choose 2}$ edges, properly placed. In \cite{BS2}, we have characterized the graphs which are periodically minimally rigid, when generically placed.

\section{Expansion and pointedness}\label{analog}

\vspace{5pt}
We begin with a lemma showing that a non-trivial (infinitesimal) deformation cannot be expansive unless all edges emanating from a vertex are {\em pointed},  for all vertices which are not in some rigid component of the framework.

\begin{lemma}\label{pointed}
Suppose the origin $0\in R^d$ is connected by edges (of fixed length) to vertices $v_i$, $i=1,\cdots,n$. If the vectors $v_i$ allow a linear combination

\begin{equation}\label{dependence}
\sum_{i=1}^n a_i v_i =0, \ \ \mbox{with all}\ \ a_i > 0 \ ,
\end{equation}

\noindent
then no infinitesimal deformation with velocity vectors $0$ at $0$ and $\dot{v}_i$ at $v_i$
can be effectively expansive. That is, the system of equalities and inequalities

\begin{equation}\label{bars}
\langle v_i,\dot{v}_i \rangle =0,\ \ i=1,\cdots,n
\end{equation}

\begin{equation}\label{expand}
\langle v_i-v_j,\dot{v}_i -\dot{v}_j \rangle \geq 0,\ \ i\neq j
\end{equation}

\noindent
with a least one strict inequality in (\ref{expand}), has no solution $\dot{v}_i$, $i=1,\cdots,n$.
\end{lemma}

\medskip \noindent
{\em Proof:}\ The infinitesimal expansion relations (\ref{expand}) imply

\begin{equation}\label{expandBis}
0 \geq \langle v_i,\dot{v}_i \rangle +\langle v_j,\dot{v}_j \rangle \geq
\langle v_i , \dot{v}_j \rangle + \langle v_j , \dot{v}_i \rangle,\ \ i\neq j
\end{equation}

\noindent
This gives

\begin{equation}\label{raa}
0 = \langle \sum_{i=1}^n a_iv_i, \sum_{i=1}^n a_i \dot{v}_i \rangle = 
\sum_{i=1}^n a_i^2 \langle v_i,\dot{v}_i \rangle + 
 \sum_{i < j} a_ia_j [ \langle v_i , \dot{v}_j \rangle + \langle v_j , \dot{v}_i \rangle ]   \ ,
\end{equation}

\noindent
but the right hand side would be strictly negative as soon as there is one strict inequality in 
(\ref{expand}).

\begin{lemma}\label{pointedBis}
Suppose the origin $0\in R^d$ is connected by edges (of fixed length) to vertices $v_i$, $i=1,\cdots,n$. Suppose that this articulated system has an effective expansive infinitesimal
deformation. Then, the convex hull of the rays $\{ av_i\ |\ a\geq 0 \}$, that is, the cone

\begin{equation}\label{cone}
{\cal  C} (v) = \{ \sum_{i=1}^n a_i v_i \ |\ a_i\geq 0 \}
\end{equation}

\noindent
has the following structure. The maximal dimension of a linear subspace contained in ${\cal C}(v)$ is $d-2$ and there are hyperplanes containing this maximal  linear part, with all remaining vectors strictly on one side. Moreover, the infinitesimal deformation is trivial on the linear part of the cone. 
\end{lemma}

\medskip \noindent
{\em Proof:}\ We apply Lemma~\ref{pointed}  to the maximal subset of vectors which allow
a zero linear combination with {\em strictly positive} coefficients. We see that the infinitesimal deformation must be trivial on this subset of vectors. Since they
span the maximal linear part of the cone ${\cal C}(v)$, what remains to be shown is that this
linear part of the cone cannot have dimension $d-1$.  

\medskip
To prove this last point by contradiction, we may suppose (by convenient labeling) that $v_1,\cdots,v_d$ give a $(d-1)$ dimensional simplex (with 0 in the interior), in the linear part of the cone and are fixed.
Then, with any other vertex $v$, the system

$$ \langle v,\dot{v} \rangle=0 $$

$$ \langle v-v_i,\dot{v} \rangle = -\langle v_i,\dot{v} \rangle  \geq 0 $$

\noindent
has only the trivial solution $\dot{v}=0$. Indeed, we have a linear combination

$$ \sum_{i=1}^d b_iv_i=0,\ \mbox{with all} \ b_i > 0,\ \ \mbox{hence}\   $$

$$  \sum_{i=1}^d b_i \langle v_i, \dot{v} \rangle =0  $$

\noindent
and the inequalities of the system must all be equalities. It follows that $\dot{v}=0$. This 
shows that an effective expansive infinitesimal deformation is incompatible with a linear 
part of dimension $d-1$ in the cone ${\cal C}(v)$.

\medskip \noindent
{\bf Remark.}\  The argument shows that expansion can result only from expansion in the
{\em orthogonal complement} of the linear part of the cone. This orthogonal complement has dimension at least two and intersects the full cone ${\cal C}(v)$ in a {\em pointed cone} i.e.
with linear part reduced to 0. This motivates the following notion.

\begin{df}\label{pointedC2}
A convex cone through 0 in $R^d$ is called pointed in codimension two when the linear part of
the cone has dimension at most $d-2$. A set of vectors with their common origin at 0 is called pointed in codimension two when contained in a convex cone pointed in codimension two.
\end{df}

\medskip \noindent
Thus, expansiveness requires pointedness in codimension two, more precisely:

\begin{theorem} If a bar-and-joint framework allows an expansive infinitesimal deformation,
then the edge vectors emanating from any vertex where the deformation is effective
must be pointed in codimension two.
\end{theorem}

\section{Periodic frameworks with expansive motions}\label{type}

\begin{figure}
 \centering
 {\includegraphics[width=1.0\textwidth]{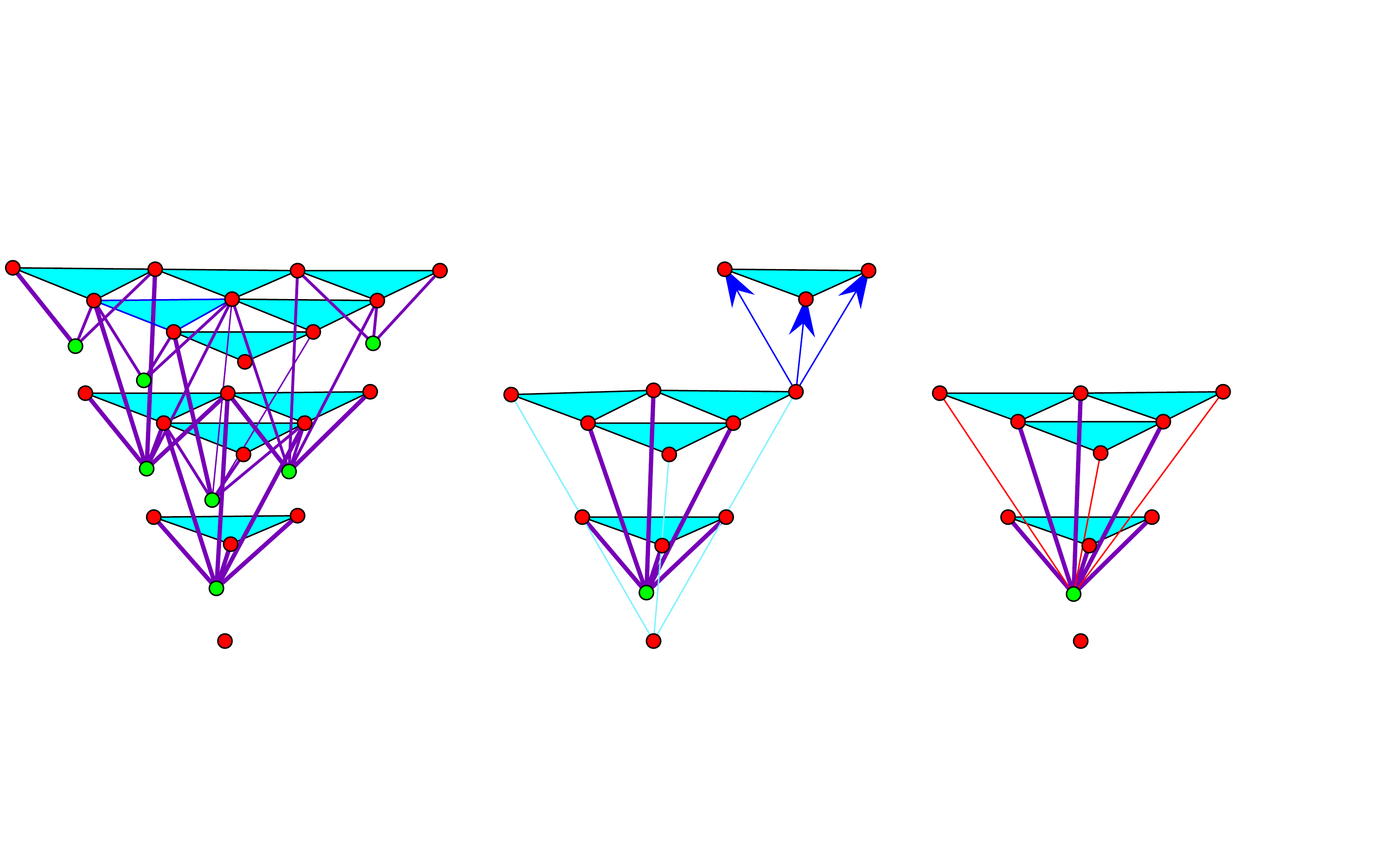}}
 \caption{ A fragment of a three dimensional periodic framework with $n=2$ vertex orbits and $m=6$ edge orbits. Essentials are isolated on the right: six edges meet at any vertex. The structure becomes periodically rigid by insertion of three new edge orbits (in red).
Blue lines and areas depict the lattice of periods on one orbit of vertices. In arbitrary dimension $d$, the analogous construction has $n=2$ and $m=d+{d\choose 2}$. This type of periodic framework has expansive capabilities in all dimensions. }
 \label{FigSimplex3D}
\end{figure}

We examine in this section a type of periodic framework which can be constructed in arbitrary
dimension $d$. The case $d=3$ is illustrated in Figure~\ref{FigSimplex3D}. For simplicity, we
may assume that we start with a {\em regular simplex} in $R^d$. The periodicity lattice will be
generated by the edge vectors of this simplex. This gives one orbit of vertices, depicted in red.
For the second orbit of vertices, depicted in green, we use a representative at the center of the simplex. The framework will have $d+{d\choose 2}$ edge orbits, with all edges between green
and red vertices.

Let us assume the initial simplex with one vertex placed at the origin and the remaining $d$ vertices at $\lambda_i$, $i=1,\cdots,d$. Then, $\lambda_i$ generate the lattice of periods and the green center will be at $v=\frac{1}{d+1}\sum_{i=1}^d \lambda_i$. We can describe a complete set of representatives for edge orbits by indicating all edges emanating from this green center: $d$ of them go to the red points of the facet $\lambda_1,\cdots,\lambda_d$ and ${d\choose 2}$ of them go to $\lambda_i+\lambda_j$, $i < j$, which should be seen as  midpoints between pairs of red vertices at $2\lambda_i$, $i=1,\cdots,d$.

 \begin{prop}\label{rigid}
By adding the $d$ edge orbits represented by edges between $v$ and $2\lambda_i$, $i=1,\cdots,d$,
the framework described above is turned into a minimally rigid framework. Thus, the original
framework has, locally, a smooth deformation space of dimension $d$.
\end{prop}

\medskip \noindent
{\em Proof.}\ A geometric argument will be shorter than a rigidity matrix rank verification.
When considering the isosceles triangle $v(2\lambda_i)(2 \lambda_j)$, the bars from $v$ to $2\lambda_i$, $\lambda_i+\lambda_j$ and $2\lambda_j$ determine the other elements of the
triangle. In particular, the period $\lambda_i-\lambda_j$ has determined length. It follows that
all the edges of the simplex with vertices at $v$ and $\lambda_i$, $i=1,\cdots,d$ are determined and
likewise for the simplex with vertices at $v$ and $2\lambda_i$, $i=1,\cdots,d$. Thus, the two simplices are rigid and must be aligned along the common altitude from $v$. By periodicity, this
determines the whole framework, up to isometry. The enhanced framework is therefore minimally
rigid, while the original framework has $d$ degrees of freedom, with a smooth local deformation
space.

\begin{figure}
 \centering
 {\includegraphics[width=0.65\textwidth]{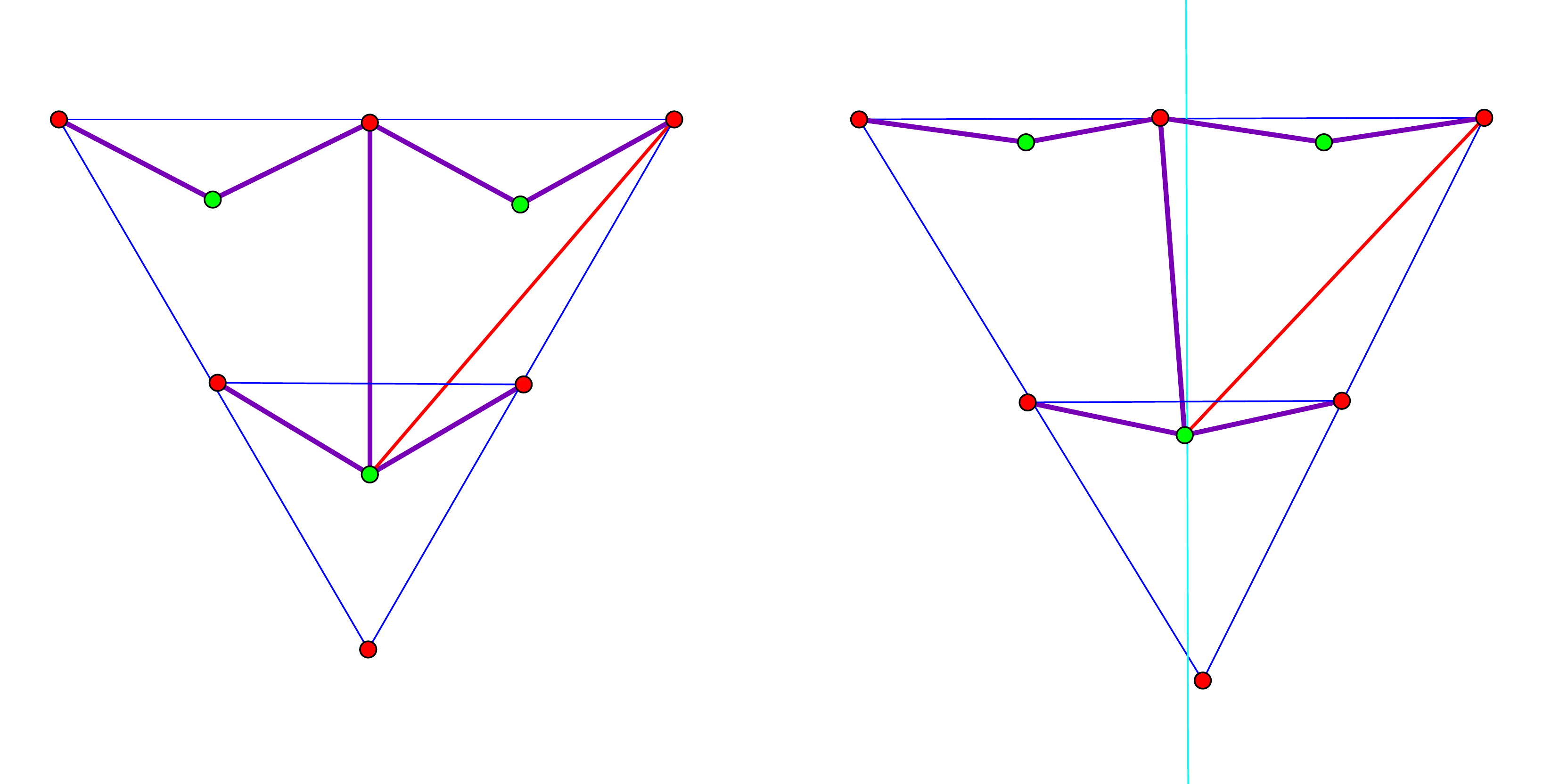}}
 \caption{Expansion illustrated for $d=2$. The absence of the edge $v(2\lambda_1)$ gives a one degree of freedom mechanism. Expansion can be observed from left to right. }
 \label{FigExp2}
\end{figure}

\begin{theorem}\label{Econe}
The minimally rigid framework in $R^d$, with edge representatives connecting $v$ with the
$2d+{d \choose 2}$ vertices at $\lambda_i, 2\lambda_i$, $i=1,\cdots,d$ and $\lambda_i+\lambda_j$, $i < j$, becomes, locally, an expansive one-degree-of-freedom
mechanism upon removal of one edge orbit of type $v(2\lambda_i)$, for any choice 
$1\leq i \leq d$. As a consequence, the framework with all  edge orbits of type $v(2\lambda_i)$, $i=1,\cdots,d$ removed, has $d$ degrees of freedom and the infinitesimal expansive cone of the framework has $d$ extremal rays, corresponding to the expansive mechanisms resulting from
single edge orbit removals.
\end{theorem}

\medskip \noindent
{\em Proof:}\ By symmetry, it is enough to examine the case of removing the edge orbit of $v(2\lambda_1)$. The planar case is shown in Figure~\ref{FigExp2}. The case of arbitrary $d$
is completely similar. The facet $(2\lambda_1)\cdots(2\lambda_d)$ expands when $(2\lambda_1)$
moves away from $(2\lambda_2)\cdots(2\lambda_d)$ along the same altitude of the facet. In order to
maintain the given bar connections, $v$ has to approach both facet 
$(2\lambda_1)\cdots(2\lambda_d)$ and its parallel $\lambda_1 \cdots\lambda_d$. An elementary
computation shows that the descent towards the latter is {\em more} than the descent towards the former, resulting in increased distance between the hyperplanes supporting these two facets.
It follows that the whole periodic motion is expansive for a sufficiently small interval of `time'.
As discussed in \cite[Borcea \& Streinu (2014b)]{BS6}, infinitesimally, expansive deformations define a convex cone in
the tangent space of the framework deformation space. For the family under consideration,
it is clear that the mechanisms corresponding to single edge orbit removals provide the extremal
rays of this cone.

\section{A stressed example}\label{stress}

\begin{figure}
 \vspace{6pt}
 {\includegraphics[width=0.95\textwidth]{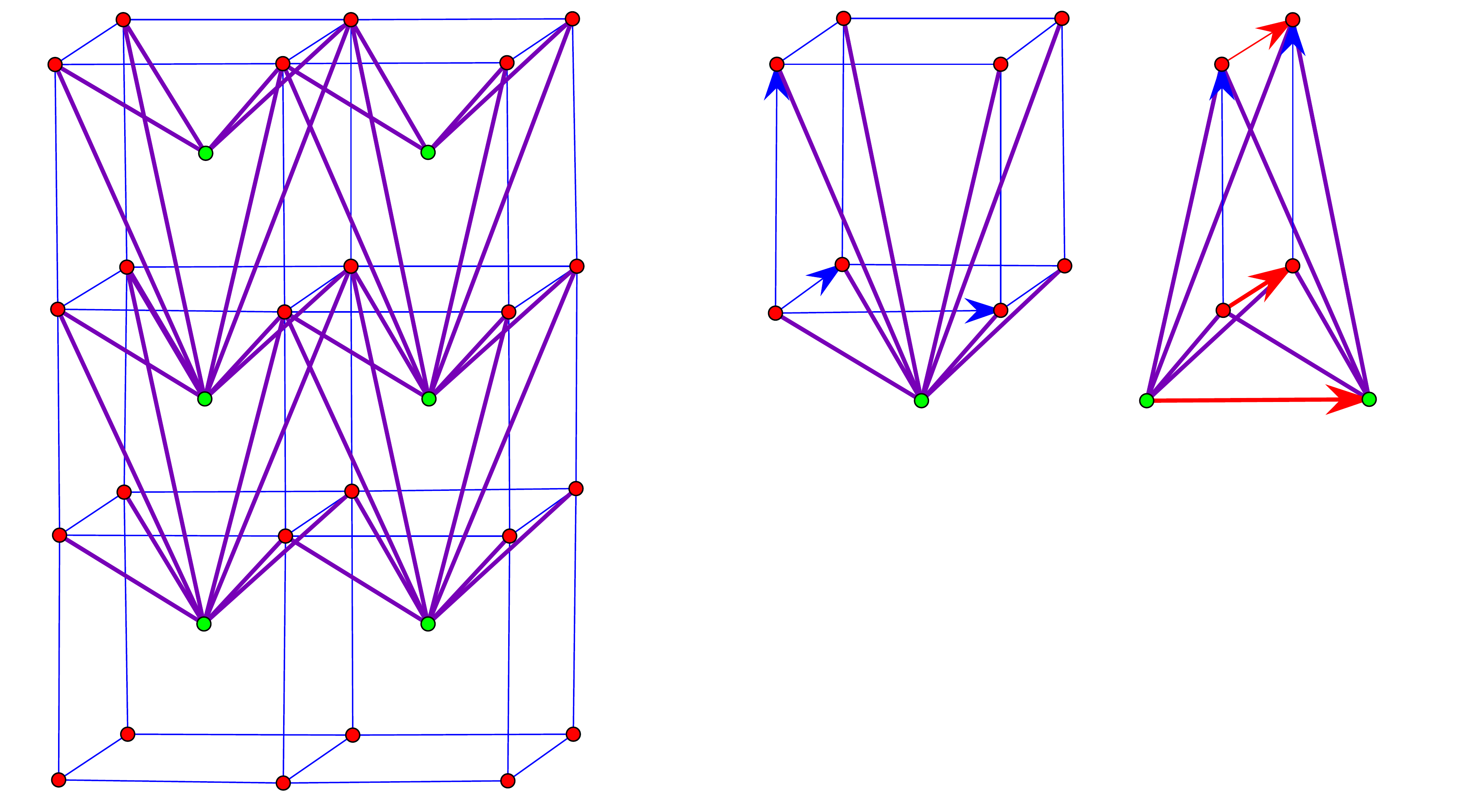}}
 \caption{ A three dimensional periodic framework with expansive capabilities. There are $n=2$ vertex orbits and $m=8$ edge orbits. Essentials are recorded in two ways on the right: eight edges  meet at any vertex.
The rightmost diagram shows that the framework can deform with two parameters (the squared lengths of the two red periods) and this signals the presence of a periodic stress. }
 \label{FigExpansiveStressed}
\end{figure}

\vspace{5pt}
In this section we consider the three dimensional periodic framework shown in Figure~\ref{FigExpansiveStressed}. In dimension two, if all rigid components are free of redundant bars, expansive frameworks do not have periodic stresses \cite[Borcea \& Streinu (2015a)]{BS4}. In higher dimensions this is not necessarily the case. Our three dimensional example will have $n=2$
vertex orbits and $m=8$ edge orbits but more than $3n+3-m=1$ degrees of freedom.
We shall identify the periodic stress which gives the framework two degrees of freedom.
Then we find the two extremal rays of the infinitesimal expansive cone.

\medskip
We place one vertex orbit representative at the origin (in red), with the lattice of periods generated by the standard basis: $\lambda_i=e_i$, $i=1,2,3$. We take a representative $v$ for the second orbit of vertices (in green) at $v=\frac{1}{2}(e_1+e_2-e_3)$. The edge orbits are represented by
the edges connecting $v$ to the following eight (red) vertices: $0, e_1, e_2, e_3$, $e_1+e_2$,
$e_2+e_3$, $e_3+e_1$, $e_1+e_2+e_3$.

\medskip \noindent
For infinitesimal considerations we may eliminate translations and assume that the origin remains fixed. Then, the unknowns of an infinitesimal deformation are: $\dot{v}$ and $\dot{\lambda}_i$, $i=1,2,3$.
The edge constraints give:

$$ \langle v, \dot{v}\rangle =0 $$

$$ \langle \lambda_i-v,\dot{\lambda}_i-\dot{v} \rangle = 0, \ \ i=1,2,3 $$

$$  \langle \lambda_i+\lambda_j-v,\dot{\lambda}_i+\dot{\lambda}_j-\dot{v} \rangle = 0, \ \ 
1\leq i < j \leq 3 $$

$$ \langle \lambda_1+\lambda_2+\lambda_3 - v, \dot{\lambda}_1+\dot{\lambda}_2+\dot{\lambda}_3-\dot{v} \rangle=0 $$

\medskip \noindent
There is a linear dependence of these conditions with coefficients 

$$ \alpha_0=-1,\ \alpha_i=1,\ \alpha_{ij}=-1, \ \alpha_{123}=1 $$

\medskip \noindent
Thus, the framework has non-trivial periodic stresses and at least two degrees of freedom,  \cite{BS1},\cite[(2015a)]{BS4}. The fact that it has exactly two is shown on the right hand side of Figure~\ref{FigExpansiveStressed}. Two local parameters which determine the periodic framework are the (squared) lengths of the periods $\lambda_1$ and $\lambda_2$. 

\begin{theorem}\label{EScone}
The three dimensional periodic framework shown in Figure~\ref{FigExpansiveStressed} is stressed and has a smooth two dimensional local deformation space. It has expansive capabilities and the two extremal rays of the infinitesimal expansive cone correspond to fixing 
the length of the period $\lambda_1$, respectively $\lambda_2$.
\end{theorem}

\medskip \noindent
We note that maintaining the length of $\lambda_1$, respectively $\lambda_2$ is equivalent
to inserting a new edge orbit represented by a connection between (red) vertices at 0 and $\lambda_1$, respectively 0 and $\lambda_2$. The expansive property of the resulting one-degree-of-freedom mechanisms is easily verified. At red vertices we'll have
aligned bars, hence not strict pointedness, but only pointedness in codimension two, as defined
above in Section~\ref{analog}.

\medskip  
{\bf Remark.}\ Other examples of stressed frameworks with expansive behavior can be obtained by relaxing the periodicity lattice of known expansive mechanisms. 

\section{Conclusions}\label{end} 

\medskip
Periodic frameworks which allow expansive motions are subject to some distinctive structural 
conditions. While complete necessary and sufficient conditions are not known for dimensions higher than two, we have identified a necessary condition of {\em pointedness in codimension two} at all vertices where expansion is effective. We have shown that expansive periodic mechanisms exist in arbitrary dimension.

\medskip 
As expected, expansive behavior is always {\em auxetic}, \cite{BS4},\cite{BS5}. Auxetic capabilities in
periodic materials and metamaterials are of considerable interest, \cite{L}, \cite{ENHR}, \cite{GGLR}, \cite{LST}. Thus, expansive structures as considered in this paper may serve, in particular, for auxetic design.

\vspace{0.7in}

Ciprian S. Borcea 

\noindent
              Department of Mathematics, Rider University, Lawrenceville, NJ 08648, USA 
              
\medskip
 Ileana Streinu 

\noindent
              Department of Computer Science, Smith College, Northampton, MA 01063, USA 
	
\end{document}